\documentclass{elsart}

\usepackage{amssymb,amsmath}

\newtheorem{theorem}{Theorem}[section]

\newtheorem{proposition}[theorem]{Proposition}

\newtheorem{question}[theorem]{Question}

\journal{Method of Func. Anal. and Top.}
\begin{document}

\begin{frontmatter}

\title{The discontinuity points set of separately continuous
functions on the products of compacts}

\author{Mykhaylyuk V.V.}

\address{Chernivtsi National University, Department of Mathematical Analysis,
Kotsjubyns'koho 2, Chernivtsi 58012, Ukraine\\vmykhaylyuk@ukr.net}

\begin{abstract}
It is solved a problem of construction of separately continuous
functions on the product of compacts with a given discontinuity
points set. We obtaine the following results.

1. For arbitrary \v{C}ech complete spaces $X$, $Y$ and a separable
compact perfect projectively nowhere dense zero set $E\subseteq
X\times Y$ there exists a separately continuous function
$f:X\times Y\to\mathbb R$ the discontinuity points set of which
equals to $E$.

2. For arbitrary \v{C}ech complete spaces $X$, $Y$ and nowhere
dense zero sets $A\subseteq X$ and $B\subseteq Y$ there exists a
separately continuous function $f:X\times Y\to\mathbb R$ such that
the projections of the discontinuity points set of $f$ coincide
with $A$ and $B$ respectively.

An example of Eberlein compacts $X$, $Y$ and nowhere dense zero
sets $A\subseteq X$ and $B\subseteq Y$ such that the discontinuity
points set of every separately continuous function $f:X\times
Y\to\mathbb R$ does not coincide with $A\times B$, and
$CH$-example of separable Valdivia compacts $X$, $Y$ and separable
nowhere dense zero sets $A\subseteq X$ and $B\subseteq Y$ such
that the discontinuity points set of every separately continuous
function $f:X\times Y\to\mathbb R$ does not coincide with $A\times
B$ are constructed.

\end{abstract}

\end{frontmatter}

\section{Introduction}
\label{Introduction}

It follows from Namioka's theorem [1] that for arbitrary compact
spaces $X$,$Y$ and a separately continuous function $f:X\times
Y\to\mathbb R$ the set $D(f)$ of discontinuity points of $f$ is a
projectively meagre set, that is $D(f)\subseteq A\times B$ where
$A\subseteq X$ and $B\subseteq Y$ are meagre sets. In this
connection, a problem on a characterization of discontinuity
points sets of separately continuous functions on the product of
two compact spaces was formulated in [2]. In other words, it is
required to establish for what projectively meagre
$F_{\sigma}$-set $E$ in the product $X\times Y$ of compact spaces
$X$ and $Y$ there exists a separately continuous function
$f:X\times Y\to\mathbb R$ with $D(f)=E$? This leads to a solving
of the inverse problem of separately continuous mappings theory of
the construction of separately continuous function with a given
discontinuity points set.

The inverse problem on $[0,1]^2$  and on the products of
metrizable spaces was studied in papers of many mathematicians
(W.~Young and G.~Young, R.~Kershner, R.~Feiock, Z.~Grande,
J.~Breckenridge and T.~Nishiura). The most general result in this
direction was obtained in [3]. It gives a characterization of
discontinuity points set for separately continuous functions of
several variables on the product of spaces each of which is the
topological product of separable metrizable factors. This result
for function of two compact variables was proved in  [4, Theorem
4] and it can be formulated in the following way.

\begin{theorem}

Let $(X_s:s\in S)$, $(Y_t:t\in T)$ be arbitrary families of
metrizable compacts, $X=\prod\limits_{s\in S}X_s$ and
$Y=\prod\limits_{t\in T}Y_t$. Then for any set $E\subseteq X\times
Y$ the following conditions are equivalent:

$(i)$\,\,\, there exists a separately continuous function
$f:X\times Y\to\mathbb R$ with $D(f)$=E;

$(ii)$\,\,\, there exists a sequence $(E_n)^{\infty}_{n=1}$ of
projectively nowhere dense zero sets $E_n\subseteq X\times Y$ such
that $E=\bigcup\limits_{n=1}^{\infty} E_n$.

\end{theorem}

Recall that a set $A$ in a topological space $X$ is called {\it a
zero set} if there exists a continuous function $f:X\to[0,1]$ such
that $A=f^{(-1)}(0)$, and {\it a co-zero set} if $A=X\setminus B$
for some zero set $B\subseteq X$. A set $E$ in the product
$X\times Y$ of topological space $X$ and $Y$ is called {\it a
projectively nowhere dense set} if $E$ is contained in the product
$A\times B$ of nowhere dense sets $A\subseteq X$ and $B\subseteq
Y$.

In other hand, the problem of construction a separately continuous
function with a given oscillation was solved in [5]. It follows
from [5] that for arbitrary separable projectively meagre
$F_{\sigma}$-set $E$ in the product $X\times Y$ of Eberlein
compacts $X$ and $Y$ there exists a separately continuous function
$f:X\times Y\to\mathbb R$ with $D(f)=E$. Besides, examples of
nonseparable closed sets $E_1$ and $E_2$ in the products of two
Eberlein compacts such that $E_1$ is the discontinuity points set
of some separately continuous function and $E_2$ is not the
discontinuity points set for every separately continuous function
on the product of the corresponding spaces.

Note that the following Price-Simon type property of Eberlein
compacts (see [6, p.170]) plays an important role in the proof of
the results of [5]. For every Eberlein compact $X$ and $x_0\in X$
there exists a sequence of nonempty open sets, which converges to
$x_0$ (a sequence $(A_n)^{\infty}_{n=1}$ of sets $A_n\subseteq Y$
{\it converges to $y_0\in Y$ in a topological space $Y$}, that is
$A_n\rightarrow x_0$, if for every neighbourhood $U$ of $y_0$ in
$Y$ there exists a number $n_0\in \mathbb N$ such that
$A_n\subseteq U$ for all $n\geq n_0$).

The problem of construction of separately continuous function on
the product of two compact spaces with a given one-point
discontinuity points set was solved in [7] using a dependence of
functions upon some quantity of coordinates technique. It was
obtained in [7] that for nonisolated points $x_0$ and $y_0$ in
compact spaces $X$ and $Y$ respectively there exists a separately
continuous function $f:X\times Y\to\mathbb R$ with
$D(f)=\{(x_0,y_0)\}$ if and only if there exist sequences
$(U_n)^{\infty}_{n=1}$ and $(V_n)^{\infty}_{n=1}$ of nonempty
co-zero sets $U_n\subseteq X$ and $V_n\subseteq Y$ which converges
to $x_0$ and $y_0$ respectively, besides, $x_0\not\in U_n$ and
$y_0\not\in V_n$ for every $n\in\mathbb N$.

Note that a solving of the inverse problem for a $F_{\sigma}$-set
$E=\bigcup\limits_{n=1}^{\infty}E_n$ is reduces to the
construction a separately continuous function $f$ with $D(f)=E_n$
where $E_n$ is a closed set. Therefore the following questions
arise naturally in a connection with results mentioned above.

\begin{question}

Let $E$ be a projectively nowhere dense zero set in the product
$X\times Y$ of compacts $X$ and $Y$. Does there exist a separately
continuous function $f:X\times Y\to\mathbb R$ with $D(f)$=E?

\end{question}

\begin{question}

Let $E$ be a projectively nowhere dense zero set in the product
$X\times Y$ of Eberlein compacts $X$ and $Y$. Does there exist a
separately continuous function $f:X\times Y\to\mathbb R$ with
$D(f)$=E?

\end{question}

\begin{question}

Let $E$ be a separable projectively nowhere dense zero set in the
product $X\times Y$ of compacts $X$ and $Y$. Does there exist a
separately continuous function $f:X\times Y\to\mathbb R$ with
$D(f)$=E?

\end{question}

Besides, theorems on characterizations of discontinuity points
sets of separately continuous functions, which were  obtained,
have been formulated in projections properties terms. Therefore it
is important to study a weak inverse problem of construction of a
separately continuous function with given projections. It is
connected with special inverse problems of construction of
separately continuous function with a given discontinuity points
set $E$ of special type ($E=A\times B$, $E=\{x_0\}\times\{y_0\}$,
etc.), which have been studied in [8, 9]. In particular, the
special inverse problem was solved in [8] in the following cases:
for a set $A\times \{y_0\}$ where $A$ is any nowhere dense zero
set in a topological space $X$ and $y_0$ is any nonisolated point
with a countable base of neighbourhoods in a completely regular
space $Y$; and for a set $A\times B$ where $A$ and $B$ are nowhere
dense zero sets in a topological space $X$ and a locally connected
space $Y$ respectively. Thus the following question arises
naturally.

\begin{question}

Let $A$, $B$ be nowhere dense zero sets in compacts $X$ and $Y$
respectively. Does there exist a separately continuous function
$f:X\times Y\to\mathbb R$ such that the projections on $X$ and $Y$
of discontinuity points set of $f$ coincide with $A$ and $B$
respectively?

\end{question}

In this paper we give positive answers to Question 1.2 if $E$ is a
separable perfect set, and to Question 1.5. Further we construct
an example which gives a negative answer to Question 1.3 (thus to
Question 1.2), and an $CH$-example which gives a negative answer
to Question 1.4.

\section{The inverse problem on the product of compacts}

Recall some definitions and introduce some notations.

A set $A$ in a topological space $X$ is called {\it a
$\overline{G}_{\delta}$-set} if there exists a sequence
$(G_n)_{n=1}^{\infty}$ of open in $X$ sets $G_n$ such that
$A=\bigcap\limits_{n=1}^{\infty}G_n$ and
$\overline{G}_{n+1}\subseteq G_n$ for every $n\in\mathbb N$ where
$\overline{B}$ means the closure of a set $B$ in the corresponding
space.

A set $A$ in a topological space $X$ is called {\it a perfect set}
if $A$ is a perfect space in the topology induced by $X$, that is
every closed in $A$ set is a $G_{\delta}$-set in $A$.

A function $f:X\to \mathbb R$ defined on a topological space $X$
is called {\it a lower semi-continuous function at an $x_0\in X$}
if for every $\varepsilon>0$ there exists a neighbourhood $U$ of
$x_0$ in $X$ such that $f(x)>f(x_0)-\varepsilon$ for any $x\in U$,
and {\it a lower semi-continuous function} if $f$ is lower
semi-continuous at any point $x\in X$.

Let $X,Y$ be arbitrary sets. The mappings ${\mbox pr}_X:X\times
Y\to X$ and ${\mbox pr}_Y:X\times Y\to Y$ are defined as follows:
${\mbox pr}_X(x,y)= x$ and ${\mbox pr}_Y(x,y)=y$ for every $x\in
X$ and $y\in Y$. Besides, let $f:X\times Y\to\mathbb R$ be a
function. For every $x_0\in X$ and $y_0\in Y$ the functions
$f^{x_0}:Y\to\mathbb R$ and $f_{y_0}:X\to\mathbb R$ are defined as
follows: $f^{x_0}(y) = f(x_0,y)$ and $f_{y_0}(x) = f(x,y_0)$ for
any $x\in X$ and $y\in Y$.

Let  $X$ be a topological space, $A\subseteq X$ and $f:X\to\mathbb
R$. The restriction of $f$ to $A$ we denote by $f|_A$. The real
$\omega_f(A)=\sup\limits_{x',x''\in A}|f(x')-f(x'')|$ is called
{\it the oscillation of $f$ on $A$}. If $x_0\in X$ and ${\mathcal
U}$ is a system of all neighborhoods of $x_0$ in $X$ then the real
$\omega_f(x_0)=\inf\limits_{U\in{\mathcal U}}\omega_f(U)$ is
called {\it the oscillation of $f$ at $x_0$}.

For a function $f:X\to\mathbb R$ defined on a set $X$ the set
${\rm supp}\,f=\{x\in X:f(x)\ne 0\}$ is called {\it a support of
$f$}.

A completely regular space $X$ is called {\it a \v{C}ech complete
space} if for every compactification $cX$ of $X$ the set $X$ is a
$G_{\delta}$-set in $cX$ (see [10, p.297).

Let $X$ be a topological space. We say that {\it a point $x_0\in
X$ has a weak Price-Simon property in $X$} if there exists a
sequence $(U_n)^{\infty}_{n=1}$ of nonempty open in $X$ sets $U_n$
such that $U_n\rightarrow x_0$, and {\it $X$ has a weak
Price-Simon property} if every point $x\in X$ has the weak
Price-Simon property in $X$.

The following result take an important place in a solving of the
inverse problem and the method of its proof is similar to the
method which was used in [11] for the product of separable
metrizable spaces.

\begin{theorem}

Let $X$, $Y$ be completely regular spaces, $A\subseteq X$ and
$B\subseteq Y$ be nowhere dense sets, $E\subseteq A\times B$ be a
$\overline G_{\delta}$-set in $Z=X\times Y$ and
$P=\{p_n:n\in\mathbb N\}\subseteq E$ be a dense in $E$ set such
that $p_n$ has the Price-Simon property in $Z$ for every
$n\in\mathbb N$. Then there exists a lower semi-continuous
separately continuous function $f:X\times Y \to Z$ such that
$D(f)=E$.

\end{theorem}

{\bf Proof.} 
Let $(G_n)^{\infty}_{n=1}$ be a sequence of open in $Z$ sets such
that $E=\bigcap\limits_{n=1}^{\infty}G_n$ and
$\overline{G}_{n+1}\subseteq G_n$ for every $n\in\mathbb N$. Since
$A$ and $B$ are nowhere dense and each point $p_n$ has the weak
Price-Simon property in $Z$, for every $n\in\mathbb N$ there exist
sequences $(U_{nk})_{k=1}^{\infty}$ and $(V_{nk})_{k=1}^{\infty}$
of nonempty open in $X$ and $Y$ sets $U_{nk}$ and $V_{nk}$
respectively such that $W_{nk}=U_{nk}\times
V_{nk}\mathop{\rightarrow}\limits_{k\to\infty} p_n$, $U_{nk}\cap A
= V_{nk}\cap B =\O$ and $W_{nk}\subseteq G_k$ for every $k\in
\mathbb N$. For every $n,k\in\mathbb N$ pick a point $z_{nk}\in
W_{nk}$ and a continuous function $f_{nk}:Z\to [0,1]$ such that
$f_{nk}(z_{nk})=1$ and $f_{nk}(z)=0$ for any $z\in Z\setminus
W_{nk}$. Show that the function $f:X\times Y\to [0,+\infty)$,
$f(x,y)=\sum\limits_{n=1}^\infty\sum\limits_{k=n}^\infty
f_{nk}(x,y)$, has the desired properties.

For every $n\in\mathbb N$ put $W_n=Z\setminus \overline{G}_n$.
Since $f_{ik}|_{W_n}=0$ for all $i\in\mathbb N$ and $k\ge n$,
$f|_{W_n}=\sum\limits_{i=1}^\infty\sum\limits_{k=i}^\infty
f_{ik}|_{W_n}=\sum\limits_{i=1}^n\sum\limits_{k=i}^n
f_{ik}|_{W_n}$. Therefore $f$ is continuous at every point of set
$\bigcup\limits_{n=1}^\infty W_n=Z\setminus E$.

Besides, since $W_{nk}\cap ((A\times Y)\cup (X\times B))=\O$ for
any $n,k\in\mathbb N$, $f^a=f_b=0$ for any $a\in A$ and $b\in B$.
Therefore, in particular, $f$ is a lower semi-continuous
separately continuous function.

It remains to show that $E\subseteq D(f)$. Since $f(p_n)=0$ for
each $n\in\mathbb N$, $f(z_{nk})\ge f_{nk}(z_{nk})=1$ for each
$k\ge n$ and $z_{nk}\mathop{\to}\limits_{k\to\infty} p_n$, $p_n\in
D(f)$, besides, $\omega_f(p_n)\ge 1$. Using the closeness of
$F=\{z\in Z:\omega_f(z)\geq 1\}$ in $Z$ and $F\subseteq D(f)$ we
obtain $E=\overline{P}\subseteq F\subseteq D(f)$.
\hfill$\diamondsuit$

For a set which is the union of a sequence of zero sets we obtain
the following solution to the inverse problem.

\begin{theorem}

Let $X$, $Y$ be completely regular spaces, $(E_n)_{n=1}^\infty$ be
a sequence of separable projectively nowhere dense
$\overline{G}_\delta$-sets $E_n$ in $X\times Y$ and
$E=\bigcup\limits_{n=1}^\infty E_n$, besides, every point of $E$
has the weak Price-Simon property in $X\times Y$. Then there
exists a separately continuous function $f:X\times Y\to\mathbb R$
such that $D(f)=E$.

\end{theorem}

{\bf Proof.} By Theorem 2.1 for every $n\in\mathbb N$ there exists a lower
semi-continuous separately continuous function $g_n:X\times
Y\to\mathbb R$ such that $D(g_n)=E_n$. Fix any strictly increasing
homeomorphism $\varphi:{\mathbb R}\to (-1,1)$. Clearly that the
functions $f_n:X\times Y\to (-1,1)$, $f_n(x,y)=\varphi(g_n(x,y))$,
are lower semi-continuous separately continuous and $D(f_n)=E_n$
for every $n\in\mathbb N$. By [12, Corollary 2.2.2] for a
separately continuous function $f:X\times Y\to\mathbb R$,
$f(x,y)=\sum\limits_{n=1}^\infty \frac{1}{2^n}f_n(x,y)$, we have
$D(f)=\bigcup\limits_{n=1}^\infty
D(f_n)=\bigcup\limits_{n=1}^\infty E_n=E$.
\hfill$\diamondsuit$

The following result gives a positive answer to Question 1.2 under
some additional conditions on $E$.

\begin{theorem}

Let $X$, $Y$ be \v{C}ech complete spaces, $(E_n)_{n=1}^\infty$ be
a sequence of separable compact perfect projectively nowhere dense
$G_\delta$-sets $E_n$ in $X\times Y$ and
$E=\bigcup\limits_{n=1}^\infty E_n$. Then there exists a
separately continuous function $f:X\times Y\to\mathbb R$ such that
$D(f)=E$.

\end{theorem}

{\bf Proof.} Let $\tilde{X}$, $\tilde{Y}$ be the Stone-\v{C}ech
compactifications of $X$ and $Y$ respectively. Since $X$ and $Y$
are \v{C}ech complete spaces, $X$ and $Y$ are $G_\delta$-sets in
$\tilde{X}$ and $\tilde{Y}$ respectively. Therefore all the sets
$E_n$ are $G_\delta$-sets in $\tilde{X}\times \tilde{Y}$. Every
one-point subset of $E_n$ is a $G_\delta$-set in the perfect
compact $E_n$. Thus every one-point subset of $E$ is a
$G_\delta$-set in the compact $\tilde{X}\times \tilde{Y}$. Hence
every point $p\in E$ has a countable base of neighbourhoods in
$\tilde{X}\times \tilde{Y}$. Then by Theorem 2.2, there exists a
separately continuous function $\tilde{f}:\tilde{X}\times
\tilde{Y}\to\mathbb R$ such that $D(\tilde{f})=E$.

Put $f=\tilde{f}|_{X\times Y}$. Clearly that $f$ is a separately
continuous function and $D(f)\subseteq E$. It remains to show that
$E\subseteq D(f)$.

Pick a point $p=(x_0,y_0)\in E$, neighbourhoods $U$ and $V$ of
$x_0$ and $y_0$ in $X$ and $Y$ respectively. Since $X$ and $Y$ are
dense in $\tilde{X}$ and $\tilde{Y}$ respectively,
$\tilde{U}=\overline{U}$ and $\tilde{V}=\overline{V}$ are
neighbourhoods of $x_0$ and $y_0$ in $\tilde{X}$ and $\tilde{Y}$
respectively. Using the separate continuity of $\tilde{f}$ we
obtain that
$\omega_{\tilde{f}}(\tilde{U}\times\tilde{V})=\omega_{\tilde{f}}(U\times
V)=\omega_{f}(U\times V)$. Therefore
$\omega_{{f}}(p)=\omega_{\tilde{f}}(p)>0$ and $p\in D(f)$.
\hfill$\diamondsuit$

Note that the \v{C}ech completeness of $X$ and $Y$ in Theorem 2.3
cannot be weaken to the complete regularity. Indeed, it was shown
in [9, Theorem 1] that an analog of this theorem for completely
regular spaces $X$, $Y$ and one-point set $E$ does not depend of
the $ZFC$-axioms.

The method which we use to solving the weak inverse problem is
similar to the method from [8]. The following proposition gives a
possibility to remove the connection type conditions.

\begin{proposition}

Let $X$ be a compact space, $A\subseteq X$ be a zero set in $X$
which is not open in $X$. Then there exists a separately
continuous function $f:X\to[0,1]$ such that $A=f^{-1}(0)$ and for
every open in $X$ set $G\supseteq A$ there exists $n_0\in\mathbb
N$ such that $\{\frac{1}{2^n}:n\geq n_0\}\subseteq f(G)$.

\end{proposition}

{\bf Proof.} Let $g:X\to[0,1]$ be a continuous function such that
$A=g^{-1}(0)$. Since $A$ is not open in $X$,
$g^{-1}([0,\varepsilon))\setminus A \ne\O$ for every
$\varepsilon>0$. Therefore there exists a sequence
$(x_n)_{n=1}^{\infty}$ of points $x_n\in X$ such that
$g(x_{n+1})<g(x_n)<1$ for every $n\in\mathbb N$ and
$\lim\limits_{n\to\infty}g(x_n)=0$. Pick any strictly increasing
continuous function $\varphi:[0,1]\to[0,1]$ such that
$\varphi(g(x_n))=\frac{1}{2^n}$. Put $f(x)=\varphi(g(x))$ for
every $x\in X$. Clearly that $f:X\to[0,1]$ is a continuous
function and $A=f^{-1}(0)$. For every $n\in\mathbb N$ put $G_n=
f^{-1}((\frac{1}{2^n},1])$. Let $G$ be an arbitrary open in $X$
set with $A\subseteq G$. Choosing a finite subcover from the open
cover $\{G\}\cup\{G_n:n\in\mathbb N\}$ of compact space $X$ we
obtain an $n_0\in\mathbb N$ such that $G\cup G_{n_0}=X$. Since
$f(x_n)=\frac{1}{2^n}\leq\frac{1}{2^{n_0}}$ for every $n\geq n_0$,
$\{\frac{1}{2^n}: n\geq n_0\} = \{f(x_n): n\geq n_0\} \subseteq
G$.
\hfill$\diamondsuit$

The following theorem gives a positive answer to Question 1.5.

\begin{theorem}

Let $X$, $Y$ be \v{C}ech complete spaces, $(A_n)_{n=1}^{\infty}$,
$(B_n)_{n=1}^{\infty}$ be sequences of nowhere dense compact
$G_{\delta}$-sets $A_n$ and $B_n$ in $X$ and $Y$ respectively,
$A=\bigcup\limits_{n=1}^{\infty}A_n$ and
$B=\bigcup\limits_{n=1}^{\infty}B_n$. Then there exists a
separately continuous function $f:X\times Y\to \mathbb R$ such
that ${\rm pr}_XD(f)=A$ and ${\rm pr}_YD(f)=B$.

\end{theorem}

{\bf Proof.} Note that it is sufficiently to prove this theorem for
nowhere dense compact $G_{\delta}$-sets $A$ and $B$ in compacts
$X$ and $Y$ respectively and a lower semi-continuous separately
continuous function $f$ analogously as in the proof of Theorem
2.3.

Since $A$ and $B$ are zero sets in $X$ and $Y$ respectively, by
Proposition 2.4 there exist continuous functions $g:X\to[0,1]$ and
$h:Y\to[0,1]$ such that $A=g^{-1}(0)$, $B=h^{-1}(0)$ and for every
open sets $G_1\supseteq A$ and $G_2\supseteq B$ in $X$ and $Y$
respectively there exists an $n_0\in\mathbb N$ such that
$\{\frac{1}{2^n}: n\geq n_0\}\subseteq g(G_1)$ and
$\{\frac{1}{2^n}: n\geq n_0\}\subseteq h(G_2)$.

Consider the function
$$
f(x,y) = \left \{\begin{array}{rr}
 \frac{2g(x)h(y)}{g^2(x) + h^2(y)},
&
 {\rm if}\quad (x,y)\not \in A\times B;
\\
  0,
&
 {\rm if}\quad (x,y)\in A\times B.
  \end{array} \right .
$$
It is easy to see that $f$ is a lower semi-continuous separately
continuous function and $D(f)\subseteq A\times B$.

Suppose that ${\rm pr}_XD(f)\ne A$, that is, there exists an
$x_0\in A\setminus {\rm pr}_XD(f)$. Since $f$ is continuous at
every point of the compact set $\{x_0\}\times B$ and $f(x_0,y)=0$
for every $y\in B$, there exists a neighbourhood $U$ of $x_0$ in
$X$ and an open in $Y$ set $G$ such that $f(x,y)<\frac{4}{5}$ for
any $x\in U$ and $y\in G$. It follows from the choice of $h$ that
there exists an $n_0\in\mathbb N$ such that $\{\frac{1}{2^n}:
n\geq n_0\}\subseteq h(G)$. Since $A=g^{-1}(0)$ is nowhere dense,
$(g^{-1}([0,\frac{1}{2^{n_0}}))\cap U)\setminus A\ne \O$, that is
there exists an $x_1\in U$ such that $g(x_1)\in
(0,\frac{1}{2^{n_0}})$. Choose an $n\geq n_0$ and a $y_1\in G$
such that $g(x_1)\in [\frac{1}{2^{n+1}}, \frac{1}{2^n})$ and
$h(y_1)=\frac{1}{2^n}$. Then
$$
f(x_1,y_1) = \frac{2g(x_1)h(y_1)}{g^2(x_1)+h^2(y_1)}\geq \frac{2
\frac{1}{2^{n+1}}\frac{1}{2^n}}{\frac{1}{4^{n+1}} + \frac{1}{4^n}}
= \frac{4}{5},
$$
but it contradicts the choice of $U$ and $G$.

The equality ${\rm pr}_YD(f)=B$ can be obtained analogously.
\hfill$\diamondsuit$

The same reasoning as after the proof of Theorem 2.3 shows that
the \v{C}ech completeness of $X$ and $Y$ in Theorem 2.5 cannot be
weaken to the complete regularity.

\section{Separately continuous functions on the product of Eberlein
compacts}

In this section we construct an example which gives a negative
answer to Question 1.3.

Recall that a compact space $X$ which is homeomorphic to some
weakly compact subset of a Banach space is called {\it an Eberlein
compact}.  The Amir-Lindenstraus theorem [13] states that a
compact $X$ is an Eberlein compact if and only if it is
homeomorphic to some compact subset of space $c_0(T)$ ($c_0(T)$ is
the space of all functions $x:T\to \mathbb R$ such that for every
$\varepsilon>0$ the set $\{t\in T: |x(t)|\geq \varepsilon\}$ is
finite with the topology of pointwise convergence on $T$).

An idea of the corresponding space construction is closely related
to the following simple fact.

\begin{proposition}

Let $f:[0,1]^2\to\mathbb R$ be a separately continuous function.
Then there exist strictly decreasing sequences
$(a_n)^{\infty}_{n=1}$ and $(b_n)^{\infty}_{n=1}$ of reals $a_n,
b_n\in (0,1]$ such that $\lim\limits_{n\to\infty} a_n =
\lim\limits_{n\to\infty} b_n = 0$ and $|f(a_n,b_m) - f(0,0)|<
\frac{1}{\min\{n,m\}}$ for every $n,m\in\mathbb N$.

\end{proposition}

{\bf Proof.} Since $f_0$ is continuous at $0$, there exists an $a_1\in
(0,1)$ such that
$$
|f(a_1,0)-f(0,0)|<\frac{1}{2}.
$$
Using the continuity of $f^0$ and $f^{a_1}$ at $0$ we choose
$b_1\in (0,1)$ such that
$$
|f(0,b_1)-f(0,0)|<\frac{1}{2}\quad\mbox{and}\quad|f(a_1,y)-f(a_1,0)|<\frac{1}{2}
$$
for every $y\in [0,b_1]$. Further, using the continuity of $f_0$
and $f_{b_1}$ at $0$ choose an $a_2\in (0,\min\{\frac{1}{2},
a_1\})$ so that
$$
|f(a_2,0)-f(0,0)|<\frac{1}{4}\quad\mbox{and}\quad|f(x,b_1)-f(0,b_1)|<\frac{1}{2}
$$
for every $x\in [0,a_2]$. Since $f^0$ and $f^{a_2}$ are continuous
at $0$, there exists $b_2\in(0,\min\{\frac{1}{2}, b_1\})$ such
that
$$
|f(0, b_2)-f(0,0)|<\frac{1}{4}\quad\mbox{and}\quad
|f(a_2,y)-f(a_2,0)|<\frac{1}{4}
$$
for every $y\in [0,b_2]$.

Continuing this procedure to infinity we obtain strictly
decreasing sequences $(a_n)_{n=1}^\infty$ and $(b_n)_{n=1}^\infty$
of reals $a_n,b_n\in(0,1]$ such that $\lim\limits_{n\to\infty} a_n
= \lim\limits_{n\to\infty} b_n = 0$ and
$$
|f(a_n,0)-f(0,0)|<\frac{1}{2n},\qquad\qquad|f(0,b_n)-f(0,0)|<\frac{1}{2n},
$$
$$
|f(a_n,y)-f(a_n,0)|<\frac{1}{2n}\qquad\mbox{and}\qquad|f(x,b_n)-f(0,b_n)|<\frac{1}{2n}
$$
for every $y\in [0,b_n]$ and $x\in [0,a_{n+1}]$. Then for $m\ge n$
we have
$$
|f(a_n,b_m)-f(0,0)|\le |f(a_n,b_m)-f(a_n,0)|+|f(a_n,0)-f(0,0)|<
$$
$$
<\frac{1}{2n}+\frac{1}{2n}=\frac{1}{n}.
$$
And for $n>m$ we have
$$
|f(a_n,b_m)-f(0,0)|\le |f(a_n,b_m)-f(0,b_m)|+|f(0,b_m)-f(0,0)|<
$$
$$
<\frac{1}{2m}+\frac{1}{2m}=\frac{1}{m}.
$$
\hfill$\diamondsuit$

For the topological product $X=\prod\limits_{s\in S}X_s$ of a
family $(X_s:s\in S)$ of topological spaces $X_s$ and a nonempty
basic open set $U=\prod\limits_{s\in S}U_s$ put $R(U)=\{s\in S:
U_s\ne X_s\}$. Let, besides, $Y$ be a subspace of $X$. An nonempty
open in $Y$ set $V$ is called {\it a basic open set} if there
exists a basic open set $U=\varphi(V)$ in $X$ such that $V=U\cap
Y$. For any nonempty basic open set $V$ in $Y$ we put
$R(V)=R(\varphi(V))$.

The following theorem is a main result of this section.

\begin{theorem}

There exist Eberlein compacts $X$ and $Y$ and nowhere dense zero
sets $A$ and $B$ in $X$ and $Y$ respectively such that $D(f)\ne
A\times B$ for every separately continuous function $f:X\times
Y\to\mathbb R$.

\end{theorem}

{\bf Proof.} Denote the set of all strictly decreasing sequences
$s=(\alpha_n)_{n=1}^\infty$ of reals $\alpha_n\in (0,1]$ such that
$\lim\limits_{n\to\infty}\alpha_{n}=0$ by $S_0$ and $S=\{0\}\cup
S_0$. For every $s=(\alpha_n)\in S_0$ and $n\in{\mathbb
N}_0={\mathbb N}\cup\{0\}$ the function $x(s,n)\in [0,1]^S$ is
defined as follows: if $n\in {\mathbb N}$, then
$$
x(s,n)(t)=\left\{\begin{array}{lll}
  1, & \mbox{if} & t=s, \\
  0, & \mbox{if} & t\in S_0\setminus \{s\}, \\
  \alpha_n, & \mbox{if} & t=0,\\
\end{array}
\right.
$$
and
$$
x(s,0)(t)=\left\{\begin{array}{lll}
  1, & \mbox{if} & t=s, \\
  0, & \mbox{if} & t\in S\setminus\{s\}.\\
\end{array}
\right.
$$

Put $X_0=[0,1]\times \{0\}^{S_0}$, $X_s=\{x(s,n):n\in {\mathbb
N}_0\}$ for every $s\in S_0$ and $X=\bigcup\limits_{s\in S} X_s$.
Show that $X$ is a closed subspace of $Z=[0,1]\times
\{0,1\}^{S_0}$.

Since for every $x\in X$ the set $\{s\in S_0:x(s)=1\}$ has at most
one element, all functions $z\in \overline{X}$ have the same
properties. Therefore it is sufficient to prove that for every
$z\in Z\setminus X$ with $|\{s\in S_0:z(s)=1\}|\le 1$ there exists
an open neighbourhood $U$ of $z$ in $Z$ such that $U\cap X=\O$.

Pick $z_0\in Z\setminus X$ such that $|\{s\in S_0:z(s)=1\}|\le 1$.
Note that $|\{s\in S_0:z(s)=1\}|= 1$. Indeed, if $z_0(s)=0$ for
every $s\in S_0$ then $z_0\in X_0$ which contradicts the choice of
$z_0$. Pick $s=(\alpha_n)\in S_0$ such that $z_0(s)=1$. Since
$z_0\ne x(s,0)$, $z_0(0)>0$. Note $z_0\ne x(s,n)$ for every
$n\in\mathbb N$, therefore $z_0(0)\ne\alpha_n$ for every
$n\in\mathbb N$. Choose an open neighbourhood $I$ of $z_0(0)$ in
[0,1] so that $0\not\in I$ and $\alpha_n\not\in I$ for every
$n\in\mathbb N$. For the open neighbourhood $U=\{z\in Z: z(0)\in
I, z(s)=1 \}$ of $z_0$ in $Z$ we have $U\cap X=\O$.

Thus $X$ is a compact. Since the supports of all functions $x\in
X$ are finite, $X$ is an Eberlein compact by [13].

Put $A=\{x\in X: x(0)=0\}$.  Clearly that $A$  is a zero set in
$X$. Since $X_s=\overline{X_s\setminus A}$ for every $s\in S$,
$X\setminus A$ is a dense in $X$ set. Therefore $A$ is a nowhere
dense in $X$ set.

Denote $Y=X$, $B=A$ and suppose that there exists a separately
continuous function $f:X\times Y\to\mathbb R$ such that
$D(f)=A\times B$. Note that the function $\varphi:X_0\to [0,1]$,
$\varphi(x)=x(0)$, is a homeomorphism, therefore the function
$g:[0,1]^2\to\mathbb R$,
$g(u,v)=f(\varphi^{-1}(u),\varphi^{-1}(v))$, is  separately
continuous. By Proposition 3.1 there exist strictly decreasing
sequences $(u_n)_{n=1}^{\infty}$ and $(v_n)_{n=1}^{\infty}$ of
reals $u_n,v_n\in (0,1]$ such that $\lim\limits_{n\to\infty} u_n =
\lim\limits_{n\to\infty} v_n = 0$ and
$|g(u_n,v_m)-g(0,0)|<\frac{1}{\min\{n,m\}}$ for every
$n,m\in\mathbb N$.

For every $n\in\mathbb N$ put $x_n=\varphi^{-1}(u_n)$ and
$y_n=\varphi^{-1}(v_n)$. Since for every $n,m\in\mathbb N$ $f$ is
a jointly continuous function at $(x_n,y_m)$, there exist basic
open neighbourhoods $U_{nm}$ and $V_{nm}$ of $x_n$ and $y_m$ in
$X$ and $Y$ respectively such that
$|f(x,y)-f(x_n,y_m)|<\frac{1}{\min\{n,m\}}$ for every $x\in
U_{nm}$ and $y\in V_{nm}$.

Consider an at most countable set $T=\bigcup\limits_{n,m=1}^\infty
(R(U_{nm})\cup R(V_{nm}))$. Since the set of all subsequences of
some sequence has the cardinality $2^{\aleph_0}$, there exists an
increasing sequence $(n_k)_{k=1}^\infty$ of $n_k\in\mathbb N$ such
that $s=(\alpha_k)_{k=1}^\infty$, $t=(\beta_k)_{k=1}^\infty\not\in
T$, where $\alpha_k=u_{n_k}$ and $\beta_k=v_{n_k}$ for every
$k\in\mathbb N$. Note that $x_0=x(s,0)\in A$ and $y_0=x(t,0)\in
B$. Show that $f$ is continuous at $(x_0,y_0)$, which is
impossible.

Fix $\varepsilon>0$ and choose a number $k_0$ such that
$\frac{1}{n_{k_0}}<\frac{\varepsilon}{2}$. Note that the sets
$U=\{x\in X:x(s)=1,
x(0)\le\alpha_{k_0}\}=\{x(s,k):k=0\quad\mbox{or}\quad k\ge k_0\}$
and $V=\{y\in Y: y(t)=1,
y(0)\le\beta_{k_0}\}=\{x(t,k):k=0\quad\mbox{or}\quad k\ge k_0\}\}$
are neighbourhoods of $x_0$ and $y_0$ in $X$ and $Y$ respectively.
Pick $i,j\ge k_0$. It follows from $s\not\in R(U_{n_in_j})$,
$t\not\in R(V_{n_in_j})$, $x_{n_i}(0)=\alpha_i$ and
$y_{n_j}(0)=\beta_j$, that $x(s,i)\in U_{n_in_j}$ and $x(t,j)\in
V_{n_in_j}$. Therefore
$|f(x(s,i),x(t,j))-f(x_{n_i},y_{n_j})|<\frac{1}{\min\{n_i,n_j\}}$.
It follows from $f(x_{n_i},y_{n_j})=g(u_{n_i},v_{n_j})$ and the
choosing of sequences $(u_n)_{n=1}^{\infty}$ and
$(v_n)^{\infty}_{n=1}$ that
$$
|f(x(s,i),x(t,j))-g(0,0)|<\frac{2}{\min\{n_i,n_j\}}.
$$
Since $f$ is a separately continuous function,
$x_0=\lim\limits_{n\to\infty} x(s,n)$ and
$y_0=\lim\limits_{n\to\infty} x(t,n)$,
$|f(x(s,i),y_0)-g(0,0)|\le\frac{2}{n_i}$,
$|f(x_0,x(t,j))-g(0,0)|\le\frac{2}{n_j}$ and $f(x_0,y_0)=g(0,0)$.
Thus $|f(x,y)-f(x_0,y_0)|\le\frac{2}{n_{k_0}}<\varepsilon$ for
every $(x,y)\in U\times V$.
\hfill$\diamondsuit$

\section{Separately continuous functions on the products of
separable Valdivia compacts}

Recall that a compact space $X$ is called {\it a Corson compact}
if it is homeomorphic to a compact $Z\subseteq {\mathbb R}^T$ such
that $|{\rm supp}\,z|\leq \aleph_0$ for every $z\in Z$, and {\it a
Valdivia compact} if it is homeomorphic to a compact $Z\subseteq
{\mathbb R}^T$ such that the set $\{z\in Z:|{\rm supp}\,z|\leq
\aleph_0\}$ is dense in $Z$. Clearly that any Corson compact is a
Valdivia compact. Besides, it follows from [13] that any Eberlein
compact is a Corson compact.

Since every separable subset of a Corson compact is metrizable, it
follows from Theorem 2.3 that Question 1.4 has a positive answer
for Corson compacts. Therefore it is naturally to establish
whether is it true for Valdivia compacts.

In this section we show that in $CH$-assumption Question 1.4 has a
negative answer even for separable Valdivia compacts.

The following notation is an important tool for the construction
of the corresponding example.

Let $X\subseteq [0,1]^S$, $Y\subseteq [0,1]^T$ be arbitrary
spaces, $s_0\in S$, $t_0\in T$ and $f:X\times Y\to \mathbb R$ be a
function. We say that sequences $(u_n)^{\infty}_{n=1}$ and
$(v_n)^{\infty}_{n=1}$ of reals $u_n, v_n \in (0,1]$ {\it nullify
$f$ in the coordinates $s_0$ and $t_0$} if the following
conditions hold:

\hangindent=1cm\noindent$(1_n)$\centerline{$|f(x,y_1)-f(x,y_2)| <
\frac{1}{n}$}\\ for every $n\in \mathbb N$, $x\in X$ with
$x(s_0)=u_n$, $m\geq n$, $y_1,y_2\in Y$ with $y_1(t)=y_2(t)$ for
$t\in T\setminus\{t_0\}$, $y_1(t_0)=v_m$ and $y_2(t_0)=0$;

\hangindent=1cm\noindent$(2_m)$\centerline{$|f(x_1,y)-f(x_2,y)| <
\frac{1}{m}$}\\ for every $m\in \mathbb N$, $y\in Y$ with
$y(t_0)=v_m$, $n>m$, $x_1,x_2\in X$ with $x_1(s)=x_2(s)$ for $s\in
S\setminus\{s_0\}$, $x_1(s_0)=u_n$ and $x_2(s_0)=0$.

\begin{proposition}

Let $X\subseteq [0,1]^S$, $Y\subseteq [0,1]^T$ be compacts,
$s_0\in S$, $t_0\in T$, $A=\{x\in X: x(s_0)=0\}$, $B=\{y\in Y:
y(t_0)=0\}$, $(a_n)^{\infty}_{n=1}$ and $(b_n)^{\infty}_{n=1}$ be
strictly decreasing sequences of reals $a_n, b_n \in (0,1]$ with
$\lim\limits_{n\to\infty} a_n = \lim\limits_{n\to\infty} b_n = 0$
and $f:X\times Y\to \mathbb R$ be a function with $D(f)\subseteq
A\times B$. Then there exist subsequences $(u_n)^{\infty}_{n=1}$
and $(v_n)^{\infty}_{n=1}$ of sequences $(a_n)^{\infty}_{n=1}$ and
$(b_n)^{\infty}_{n=1}$ respectively which nullify $f$ in the
coordinates $s_0$ and $t_0$.

\end{proposition}

{\bf Proof.} For every $n\in\mathbb N$ put $A_n=\{x\in X: x(s_0)=a_n\}$
and $B_n=\{y\in Y: y(t_0)=b_n\}$. Since $f$ is jointly continuous
at any point of compacts $A_n\times Y$, for every $n\in\mathbb N$
there exists $\varepsilon_n>0$ such that
$$
|f(x,y_1)-f(x,y_2)| < \frac{1}{n}
$$
for any $x\in A_n$, $y_1,y_2\in Y$ with $y_1(t)=y_2(t)$ for $t\in
T\setminus\{t_0\}$ and $|y_1(t_0) - y_2(t_0)|<\varepsilon_n$.

Analogously, for every $m\in\mathbb N$ the joint continuity of $f$
on the compact $X\times B_m$ implies the existence of a
$\delta_m>0$ such that
$$
|f(x_1,y)-f(x_2,y)| < \frac{1}{m}
$$
for any $y\in B_m$, $x_1, x_2 \in X$ with $x_1(s)=x_2(s)$ for
$s\in S\setminus \{s_0\}$ and $|x_1(s_0) - x_2(s_0)|<\delta_m$.

Denote $i_1=1$ and choose strictly increasing sequences
$(i_n)_{n=2}^{\infty}$ and $(j_n)_{n=1}^{\infty}$ of numbers
$i_n,j_n \in \mathbb N$ such that $b_{j_n}<\varepsilon_{i_n}$ and
$a_{i_{n+1}}<\delta_{j_n}$ for every $n\in\mathbb N$.

It remains to put $u_n=a_{i_n}$ and $v_n=b_{j_n}$ for every
$n\in\mathbb N$.
\hfill$\diamondsuit$

Now describe a method of the construction of Valdivia compacts.

Let ${\mathcal A}$ be a system of sets $A\subseteq [0,1]$,
$S=\{0\}\cup{\mathcal A}$, $X_0=[0,1]$, $X_s=\{0,1\}$ for every
$s\in {\mathcal A}$ and $X=\prod\limits_{s\in S}X_s$. For every
finite set $T\subseteq{\mathcal A}$ put
$$
Z_T=\{x\in X: x(s)=1\,\forall s\in T, \,
x(s)=0\,\forall\,s\in{\mathcal A}\setminus
T,\,x(0)\in\bigcap\limits_{A\in T}A\},
$$
if $T\ne\O$, and
$$
Z_{\O}=\{x\in X: x(s)=0\,\forall s\in {\mathcal
A},\,x(0)\in\bigcup{\mathcal A}\}.
$$
The compact subspace $X_{\mathcal A}=\overline{\cup\{Z_T:
\,T\subseteq{\mathcal A},\,T \mbox{\,\,is\,\,finite}\}}$ of the
space $X$ is called {\it a compact generated by the system
${\mathcal A}$}. Clearly that $X_{\mathcal A}$ is a Valdivia
compact.

We use the following properties of compacts generated by systems.

\begin{proposition}

Let ${\mathcal A}$ be a system of sets $A\subseteq [0,1)$,
$X=X_{\mathcal A}$ and $s_0=A_0\in {\mathcal A}$. Then for every
$x\in X_{\mathcal A}$ if $x(s_0)=1$ then $x(0)\in \overline{A_0}$.

\end{proposition}

{\bf Proof.} It follows from the definition of $X_{\mathcal A}$ that for
every $x\in \cup\{Z_T: \,T\subseteq{\mathcal A},\,T
\mbox{\,\,is\,\,finite}\}$ if $x(s_0)=1$ then $x(0)\in A_0$. It
remains to apply the closure operation.
\hfill$\diamondsuit$

\begin{proposition}

Let ${\mathcal A}$ be a system of sets $A\subseteq [0,1)$ such
that the set $A_0 =\bigcup{\mathcal A}$ is at most countable. Then
$X_{\mathcal A}$ is a separable compact.

\end{proposition}

{\bf Proof.} Let $A_0=\{a_n:n\in\mathbb N\}$. For every $n\in\mathbb N$
put ${\mathcal A}_n=\{A\in{\mathcal A}:a_n\in A\}$ and $X_n=\{x\in
X_{\mathcal A}: x(0)=a_n,\, x(s)=0\,\forall s\in {\mathcal
A}\setminus{\mathcal A}_n\}$. Note that for every finite set
$T\subseteq {\mathcal A}_n$ the function
$$
x_T(s)= \left \{\begin{array}{rr}
 a_n,
&
 {\rm if}\quad s=0;
\\
  1,
&
 {\rm if}\quad s\in T;
\\
  0,
&
 {\rm if}\quad s\in {\mathcal A}\setminus T,
  \end{array} \right .
$$
belongs to $X_{\mathcal A}$. Therefore
$X_n=\{a_n\}\times\prod\limits_{s\in {\mathcal
A}_n}\{0,1\}\times\prod\limits_{s\in {\mathcal
A}\setminus{\mathcal A}_n}\{0\}$. Since $|{\mathcal A}_n|\leq
2^{\aleph_0}$, $X_n$ is a separable space by
Hewitt-Marczewski-Pondiczery theorem [10, p.133].

Since $Z_T\subseteq\bigcup\limits_{n=1}^{\infty}X_n$ for any
finite set $T\subseteq {\mathcal A}$, $X_{\mathcal
A}=\overline{\bigcup\limits_{n=1}^{\infty}X_n}$. Thus $X_{\mathcal
A}$ is a separable space.
\hfill$\diamondsuit$

\begin{proposition}

Let ${\mathcal A}$ be a system of sets $A\subseteq [0,1)$ and
${\mathcal B}\subseteq {\mathcal A}$ such that $\bigcup{\mathcal
B} = \bigcup{\mathcal A}$. Then $\varphi(X_{\mathcal A})=
X_{\mathcal B}$ where $\varphi:X_{\mathcal A}\to \mathbb
R^{\{0\}\cup{\mathcal B}}$, $\varphi(x)=x|_{\{0\}\cup{\mathcal
B}}$.

\end{proposition}

{\bf Proof.} The inclusion $X_{\mathcal B}\subseteq \varphi(X_{\mathcal
A})$ follows immediately from the definition of a compact
generated by a system.

Fix a finite subsystem $T\subseteq{\mathcal A}$. If
$T\cap{\mathcal B}=\O$, then $\bigcup{\mathcal B} =
\bigcup{\mathcal A}$ implies $\varphi(Z_T)\subseteq X_{\mathcal
B}$. If $T\cap{\mathcal B}\ne\O$, then the inclusion
$\varphi(Z_T)\subseteq X_{\mathcal B}$ follows from the definition
of $Z_T$.

Since the set $\bigcup\{\varphi(Z_T): \,T\subseteq{\mathcal A},\,T
\mbox{\,\,is \,\,finite}\}$ is dense in $\varphi(X_{\mathcal A})$,
$X_{\mathcal B}$ is dense in $\varphi(X_{\mathcal A})$. Thus
$\varphi(X_{\mathcal A})\subseteq X_{\mathcal B}$.
\hfill$\diamondsuit$

We use also the following two facts.

\begin{proposition}

Let $(A_n)_{n=1}^{\infty}$ be a sequence of infinite sets
$A_n\subseteq \mathbb N$ such that $\bigcap\limits_{k=1}^{n} A_k$
is an infinite set for every $n\in\mathbb N$. Then there exists an
infinite set $B\subseteq \mathbb N$ such that $|B\setminus
A_n|<\aleph_0$ for every $n\in\mathbb N$.

\end{proposition}

{\bf Proof.} It is sufficient to put $n_1=\min A_1$, $n_k =
\min(\bigcap\limits_{i=1}^{k+1}A_i\setminus \{n_1, \dots
,n_{k-1}\})$ for every $k\geq 2$ and $B=\{n_k: k\in\mathbb N\}$.
\hfill$\diamondsuit$

\begin{proposition}

Let $\omega$ be the first ordinal of some infinite cardinality.
Then there exists a bijection $\varphi:[1,\omega)^2\to[1,\omega)$
such that $\varphi(\xi,\eta)\geq\xi$ for every
$\xi,\eta\in[1,\omega)$.

\end{proposition}

{\bf Proof.} Note that $|[1,\omega)^2|=|[1,\omega)|$ by
Hessenberg's theorem [14, p.~284], that is, there exists a
bijection $\psi:[1,\omega)\to [1,\omega)^2$. For every $\xi\in
[1,\omega)$ denote $(\alpha_{\xi},\beta_{\xi})=\psi(\xi)$.

Using the transfinite induction we construct a bijection
$\tilde{\varphi}:[1,\omega)\to [1,\omega)$ such that
$\tilde{\varphi}(\xi)\geq \alpha_{\xi}$ for every $\xi\in
[1,\omega)$.

Put $\tilde{\varphi}(1)=\alpha_1$.

Assume that $\tilde{\varphi}(\eta)$ is defined for all $\eta\in
[1,\xi)$ where $\xi\in (1,\omega)$. Put $\tilde{\varphi}(\xi)
=\min([\alpha_{\xi},\omega)\setminus
\{\tilde{\varphi}(\eta):1\leq\eta<\xi\})$.

Clearly that $\tilde{\varphi}$ is an injection and
$\tilde{\varphi}(\xi)\geq \alpha_{\xi}$ for every $\xi\in
[1,\omega)$. Show that $\tilde{\varphi}$ is a surjection.

Fix a $\xi\in [1,\omega)$. Choose an $\eta\in[1,\omega)$ such that
$\psi(\eta)=(\xi,1)$, that is $a_{\eta}=\xi$ and $b_{\eta}=1$. If
$\tilde{\varphi}(\eta)\ne \xi$, then $\xi\in
\{\tilde{\varphi}(\zeta): 1\leq \zeta \leq \eta\}$.

It remains to put $\varphi = \tilde{\varphi}\circ \psi^{-1}$.
\hfill$\diamondsuit$

Let $Z$, $S$ be arbitrary sets, $X\subseteq {\mathbb R}^S$ and
$f:X\to Z$. We say that {\it $f$ depends upon a countable quantity
of coordinates } if there exists  an at most countable set
$T\subseteq S$ such that $f(x')=f(x'')$ for every $x', x''\in X$
with $x'|_T = x''|_T$. It is easy to see that for any compact
$X\subseteq {\mathbb R}^S$ every continuous function
$f:X\to\mathbb R$ depends upon a countable quantity of
coordinates. It follows from [7, Theorem 1] that if $X\subseteq
{\mathbb R}^S$ and $Y\subseteq {\mathbb R}^T$ are separable
compacts, then every separately continuous function $f:X\times
Y\to\mathbb R$ depends upon countable quantity of coordinates as a
mapping defined on $X\times Y$, that is there exist at most
countable sets $S_0\subseteq S$ and $T_0\subseteq T$ such that
$f(x',y')=f(x'',y'')$ for every $x', x''\in X$ with $x'|_{S_0} =
x''|_{S_0}$ and $y', y''\in Y$ with $y'|_{T_0} = y''|_{T_0}$.

The following theorem is the main result of this section.

\begin{theorem} {\bf (CH)}  There exist separable Valdivia
compacts $X$ and $Y$, nowhere dense separable zero sets $E$ and
$F$ in $X$ and $Y$ respectively such that $D(f)\ne E\times F$ for
every separately continuous function $f:X\times Y\to\mathbb R$.
\end{theorem}

{\bf Proof.} Put $A_0=B_0=\{0\}\cup\{\frac{1}{n}: n\in\mathbb N\}$. Using
the transfinite induction we construct families $(A_{\xi}:
1\leq\xi<\omega_1)$ and $(B_{\xi}: 1\leq\xi<\omega_1)$ of sets
$A_{\xi}=\{0\}\cup\{a^{(\xi)}_n:n\in\mathbb N\}\subseteq A_0$ and
$B_{\xi}=\{0\}\cup\{b^{(\xi)}_n:n\in\mathbb N\}\subseteq B_0$
where $(a^{(\xi)}_n)^{\infty}_{n=1}$ and
$(b^{(\xi)}_n)^{\infty}_{n=1}$ are strictly decreasing sequences
which satisfy the following conditions:

$(1)$ \,\,\,$A_{\xi}\setminus A_{\eta}$ and $B_{\xi}\setminus
B_{\eta}$ are finite sets for every $0\leq\eta <\xi<\omega_1$;

$(2)$\,\,\,for every $\xi\in [1,\omega_1)$ and separately
continuous function $g:X_{{\mathcal A}_\xi}\times X_{{\mathcal
B}_\xi}\to\mathbb R$ with $D(g)\subseteq E_{\xi}\times F_{\xi}$,
where ${\mathcal A}_\xi =\{A_{\zeta}: 0\leq\zeta<\xi\}$,
${\mathcal B}_\xi =\{B_{\zeta}: 0\leq\zeta<\xi\}$, $E_{\xi}=\{x\in
X_{{\mathcal A}_\xi}: x(0)=0\}$ and $F_{\xi}=\{y\in X_{{\mathcal
B}_\xi}: y(0)=0\}$, there exists an $\eta\in [1,\omega_1)$ such
that the sequences $(a^{(\eta)}_n)^{\infty}_{n=1}$ and
$(b^{(\eta)}_n)^{\infty}_{n=1}$ nullify $g$ in the coordinates
$s_0=0$ and $t_0=0$.

Using Proposition 4.6 choose a bijection
$$
[1,\omega_1)\ni\xi\stackrel{\varphi}\mapsto(\varphi_1(\xi),\varphi_2(\xi))
\in[1,\omega_1)^2
$$
so that $\varphi_1(\xi)\leq\xi$ for every $\xi\in[1,\omega_1)$, in
particular, $\varphi_1(1)=1$.

Since $X_{{\mathcal A}_1}$ and $X_{{\mathcal B}_1}$ are separable
by Proposition 4.3, every separately continuous function
$g:X_{{\mathcal A}_1}\times X_{{\mathcal B}_1}\to\mathbb R$ is
determined by its values on some at most countable dense subset of
$X_{{\mathcal A}_1}\times X_{{\mathcal B}_1}$. Therefore the
system ${\mathcal F}_1$ of all separately continuous functions
$g:X_{{\mathcal A}_1}\times X_{{\mathcal B}_1}\to\mathbb R$ with
$D(g)\subseteq E_1\times F_1$ has the cardinality $2^{\aleph_0}$,
that is ${\mathcal F}_1=\{g_{(1,\eta)}: 1\leq\eta<\omega_1\}$.
Using Proposition 4.1 choose subsequences
$(a_n^{(1)})^{\infty}_{n=1}$ and $(b_n^{(1)})^{\infty}_{n=1}$ of
the sequence $(\frac{1}{n})^{\infty}_{n=1}$ which nullify
$g_{\varphi(1)}$ in the coordinates $s_0=0$ and $t_0=0$.

Assume that the sets $A_{\eta}$ and $B_{\eta}$ for
$1\leq\eta<\xi<\omega_1$ are constructed such that condition $(1)$
holds and for every $\eta\in [1,\xi)$ the sequences
$(a_n^{(\eta)})^{\infty}_{n=1}$ and
$(b_n^{(\eta)})^{\infty}_{n=1}$ nullify $g_{\varphi(\eta)}$ in the
coordinates $s_0=0$ and $t_0=0$ where ${\mathcal
F}_{\eta}=\{g_{(\eta,\zeta)}: 1\leq\zeta<\omega_1\}$ is the system
of all separately continuous functions $g:X_{{\mathcal
A}_{\eta}}\times X_{{\mathcal B}_{\eta}}\to\mathbb R$ with
$D(g)\subseteq E_{\eta}\times F_{\eta}$.

It follows from Proposition 4.3 that $X_{{\mathcal A}_\xi}$ and $
X_{{\mathcal B}_\xi}$ are separable. Therefore the system
${\mathcal F}_{\xi}$ of all separately continuous functions
$g:X_{{\mathcal A}_{\xi}}\times X_{{\mathcal B}_{\xi}}\to\mathbb
R$ with $D(g)\subseteq E_{\xi}\times F_{\xi}$ has the cardinality
$2^{\aleph_0}$, that is ${\mathcal F}_{\xi}=\{g_{(\xi,\eta)}:
1\leq\eta<\omega_1\}$. Besides, since $\varphi_1(\xi)\leq\xi$, we
have $g_{\varphi(\xi)}\in \bigcup\limits^{\xi}_{\eta=1}{\mathcal
F}_{\eta}$. Using Proposition 4.5 choose strictly decreasing
sequences $(a_n)^{\infty}_{n=1}$ and $(b_n)^{\infty}_{n=1}$ of
reals $a_n, b_n\in A_0$ such that for every $\eta\in[1,\xi)$ the
sets $\{a_n:n\in\mathbb N\}\setminus A_{\eta}$ and
$\{b_n:n\in\mathbb N\}\setminus B_{\eta}$ are finite. Now using
Proposition 4.1 choose subsequences $(a_n^{(\xi)})^{\infty}_{n=1}$
and $(b_n^{(\xi)})^{\infty}_{n=1}$ of sequences
$(a_n)^{\infty}_{n=1}$ and $(b_n)^{\infty}_{n=1}$ respectively
which nullify $g_{\varphi(\xi)}$ in the coordinates $s_0=0$ and
$t_0=0$.

Clearly that the families $(A_{\xi}: 1\leq\xi<\omega_1)$ and
$(B_{\xi}: 1\leq\xi<\omega_1)$ satisfy $(1)$. Show that the
condition $(2)$ holds.

Fix  $\xi\in [1,\omega_1)$ and $g\in {\mathcal F}_{\xi}$. Then
there exists a $\zeta\in[1,\omega_1)$ such that
$g=g_{(\xi,\zeta)}$. Since $\varphi$ is a bijection, there exists
$\eta\in[1,\omega_1)$ such that $\varphi(\eta)=(\xi,\zeta)$. Then
the sequences $(a_n^{(\eta)})^{\infty}_{n=1}$ and
$(b_n^{(\eta)})^{\infty}_{n=1}$ nullify $g_{(\xi,\zeta)}$ in the
coordinates $s_0=0$ and $t_0=0$.

Put ${\mathcal A}=\{A_{\xi}:0\leq\xi<\omega_1\}$, ${\mathcal
B}=\{B_{\xi}:0\leq\xi<\omega_1\}$, $X=X_{\mathcal A}$,
$Y=Y_{\mathcal B}$, $E=\{x\in X: x(0)=0\}$ and $F=\{y\in Y:
y(0)=0\}$. Note that compacts $E$ and $F$ are homeomorphic to
$\{0,1\}^{\omega_1}$. Therefore $E$ and $F$ are separable. It easy
to see that $E$ and $F$ are nowhere dense in $X$ and $Y$
respectively. Besides, by Proposition 4.4 for every $\xi\in
[1,\omega_1)$ we have $\pi^{(\xi)}_1(X)= X_{{\mathcal A}_{\xi}}$
and $\pi^{(\xi)}_2(Y)= X_{{\mathcal B}_{\xi}}$ where
$\pi^{(\xi)}_1:X\to \mathbb R^{\{0\}\cup{\mathcal A}_{\xi}}$,
$\pi^{(\xi)}_1(x)=x|_{\{0\}\cup{\mathcal A}_{\xi}}$, and
$\pi^{(\xi)}_2:Y\to \mathbb R^{\{0\}\cup{\mathcal B}_{\xi}}$,
$\pi^{(\xi)}_2(y)=y|_{\{0\}\cup{\mathcal B}_{\xi}}$.

Suppose that $f:X\times Y\to\mathbb R$ is a separately continuous
function with $D(f)= E\times F$. Since $X$ and $Y$ are separable,
$f$ depends upon a countable quantity of coordinates, that is,
there exist a $\xi\in[1,\omega_1)$ and a function $g:X_{{\mathcal
A}_{\xi}}\times X_{{\mathcal B}_{\xi}}\to\mathbb R$ such that
$f(x,y)= g(\pi^{(\xi)}_1(x),\pi^{(\xi)}_2(y))$ for every $x\in X$
and $y\in Y$. Note that mappings $\pi^{(\xi)}_1$ and
$\pi^{(\xi)}_2$ are perfect, therefore $g$ is a separately
continuous function and $D(g)= E_{\xi}\times F_{\xi}$ by [7,
Proposition 2]. Thus, $g\in {\mathcal F}_{\xi}$. Using $(2)$
choose an $\eta\in[\xi,\omega_1)$ such that the sequences
$(a_n^{(\eta)})^{\infty}_{n=1}$ and
$(b_n^{(\eta)})^{\infty}_{n=1}$ nullify $g$ in the coordinates
$s_0=0$ and $t_0=0$.

Put $s_1=A_{\eta}$, $t_1=B_{\eta}$, $u_n=a_n^{(\eta)}$ and
$v_n=b_n^{(\eta)}$ for every $n\in\mathbb N$. Clearly that the
sequences $(u_n)^{\infty}_{n=1}$ and $(v_n)^{\infty}_{n=1}$
nullify $f$ in the coordinates $s_0=0$ and $t_0=0$.

It follows from Namioka's theorem [1] that for the separately
continuous functions $h_1:E\times Y\to\mathbb R$, $h_1=f|_{E\times
Y}$, and $h_2:X\times F\to\mathbb R$, $h_2=f|_{X\times F}$, there
exist dense in $E$ and $F$ respectively $G_{\delta}$-sets
$E_0\subseteq E$ and $F_0\subseteq F$ such that $h_1$ is jointly
continuous at any point of $E_0\times Y$ and $h_2$ is jointly
continuous at any point of $X\times F_0$. Note that the sets
$\{x\in E:x(s_1)=1\}$ and $\{y\in F:y(t_1)=1\}$ are open and
nonempty in $E$ and $F$ respectively. Therefore there exist an
$x_0\in E_0$ and a $y_0\in F_0$ such that $x_0(s_1)=1$ and
$y_0(t_1)=1$.

Show that $f$ is jointly continuous at $(x_0,y_0)$.

Fix $\varepsilon>0$ and $k\in \mathbb N$ so that $\frac{1}{k}\leq
\frac{\varepsilon}{2}$. Using the continuity of $h_1$ and $h_2$ at
$(x_0,y_0)$ choose $l\in\mathbb N$, $s_2,\dots,s_l\in {\mathcal
A}$, $t_2,\dots,t_l\in {\mathcal B}$ and $\delta<\min\{u_k,v_k\}$
such that
$$
|f(x,y)-f(x_0,y_0)|<\frac{\varepsilon}{2}
$$
for every $(x,y)\in ((U\cap E)\times V)\cup (U\times (V\cap F))$,
where $U=\{x\in X: x(0)<\delta, x(s_i)=x_0(s_i)\mbox{ for } 1\le
i\le l\}$ and $V=\{y\in Y: y(0)<\delta, y(t_i)=y_0(t_i)\mbox{ for
} 1\le i\le l\}$.

Show that $|f(x,y)-f(x_0,y_0)|<\varepsilon$ for every $x\in U$ and
$y\in V$. Fix $x\in U$ and $y\in V$. Clearly that it is sufficient
to consider the case of $x(0)>0$ and $y(0)>0$. Note that
Proposition~4.2 implies $x(0)\in A_\eta$ and $y(0)\in B_\eta$. It
follows from the choice of $\delta$ that there exist $n,m\ge k$
such that $x(0)=u_n$ and $y(0)=v_m$.

Assume that $m\ge n$. Since the function
$$
\tilde{y}(t)=\left\{
\begin{array}{ll}
  0, & \mbox{if } t=0, \\
  y(t), & \mbox{if } t\in {\mathcal B},\\
\end{array}\right.
$$
belongs to $V\cap F$ and the sequences $(u_i)_{i=1}^\infty$ and
$(v_i)_{i=1}^\infty$ nullify $f$ in the coordinates $s_0=0$ and
$t_0=0$, we have
$$
|f(x,y)-f(x,\tilde{y})|<\frac{1}{n}\le\frac{1}{k}\le\frac{\varepsilon}{2}.
$$
Then
$$
|f(x,y) - f(x_0,y_0)|\leq |f(x,y_0) - f(x,\tilde{y})| +
|f(x,\tilde{y}) - f(x_0,y_0)| < \frac{\varepsilon}{2} +
\frac{\varepsilon}{2}=\varepsilon.
$$

If $n>m$, then we reason analogously.

Thus $f$ is jointly continuous at $(x_0,y_0)$, which is
impossible.
\hfill$\diamondsuit$

\end{document}